\documentclass[11pt]{article}
\usepackage{amssymb}
\usepackage{epsfig}
\usepackage{amsmath}
\usepackage[margin = 1 in]{geometry}
\usepackage{bm}
\newtheorem{theorem}{Theorem}[section]

\newtheorem{lemma}[theorem]{Lemma}
\newtheorem{corollary}[theorem]{Corollary}
\newtheorem{conjecture}[theorem]{Conjecture}
\newtheorem{remark}[theorem]{Remark}

\newcommand{\reals}{\ensuremath{\mathbb{R}}}

\newcommand{\brac}[1]{\left(#1\right)}
\newcommand{\bfrac}[2]{\brac{\frac{#1}{#2}}}
\newcommand{\beq}[1]{\begin{equation}\label{#1}}
\newcommand{\eeq}{\end{equation}}
\newcommand{\blem}[1]{\begin{lemma}\label{#1}}
\newcommand{\elem}{\end{lemma}}
\newcommand{\bth}[1]{\begin{theorem}\label{#1}}
\newcommand{\enth}{\end{theorem}}
\newcommand{\brem}[1]{\begin{remark}\label{#1}}
\newcommand{\erem}{\end{remark}}

\def\cP{{\cal P}}

\def\a{\alpha}
\def\b{\beta}
\def\d{\delta}
\def\D{\Delta}
\def\e{\varepsilon}
\def\f{\phi}
\def\g{\gamma}

\def\k{\kappa}

\def\th{\theta}

\def\l{\lambda}
\def\m{\mu}
\def\n{\nu}
\def\p{\pi}

\def\r{\rho}
\def\R{\Rho}
\def\s{\sigma}
\def\S{\Sigma}
\def\t{\tau}
\def\om{\omega}
\def\OM{\Omega}

\def\cF{{\cal F}}
\def\para{\vspace{.1in}}

\def\whp{{\bf whp}}
\def\Whp{{\bf Whp}}

\newcommand{\rdown}[1]{\left\lfloor #1 \right\rfloor }

\newcommand{\set}[1]{\left\{#1\right\}}
\def\vol{{\bf vol}_N}
\def\diam{{\bf diam}}
\newcommand{\proofstart}{{\bf Proof\hspace{2em}}}
\newcommand{\proofend}{\hspace*{\fill}\mbox{$\Box$}}
\def\eps{\varepsilon}
\def\R{\reals}
\def\E{{\sf E}}
\def\Pr{{\sf P}}
\def\ba{{\bf a}}
\def\bp{{\bf p}}

\newcommand{\ignore}[1]{}
\def\sM{\s_{\max}}
\def\sm{\s_{\min}}
\def\hf{\hat{f}}
\def\cB{{\cal B}}
\def\hg{\hat{g}}

\def\cE{{\cal E}}

\begin{document}
\title{\Large Logconcave Random Graphs}
\author{Alan Frieze\thanks{Department of Mathematical Sciences, Carnegie Mellon University, Pittsburgh PA15213, email: alan@random.math.cmu.edu.
Research supported in part by NSF award CCF-0502793.}
\and Santosh Vempala\thanks{School of Computer Science, Georgia Tech, Atlanta GA 30332, email: vempala@cc.gatech.edu. Supported in part by NSF award CCF-0721503.}
\and
Juan Vera\thanks{Department of Management Sciences,
University of Waterloo, jvera@uwaterloo.ca}
}

\date{}
\maketitle

\begin{abstract}
We propose the following model of a random graph on $n$ vertices.
Let $F$ be a distribution in $R_+^{n(n-1)/2}$ with a coordinate for
every pair $ij$ with $1 \le i,j \le n$. Then $G_{F,p}$ is the
distribution on graphs with $n$ vertices obtained by picking a
random point $X$ from $F$ and defining a graph on $n$ vertices whose edges are pairs
$ij$ for which $X_{ij} \le p$. The standard Erd\H{o}s-R\'{e}nyi model is the
special case when $F$ is uniform on the $0$-$1$ unit cube.  We examine basic properties such as
the connectivity threshold for quite general distributions. We also consider cases where the $X_{ij}$
are the edge weights in some random instance of a combinatorial optimization problem.
By choosing suitable distributions, we
can capture random graphs with interesting properties such as
triangle-free random graphs and weighted random graphs with bounded
total weight.
\end{abstract}
\thispagestyle{empty} \setcounter{page}{0} \clearpage

\section{Introduction}
Probabilistic combinatorics is today a thriving field bridging the
classical area of probability with modern developments in
combinatorics. The theory of random graphs, pioneered by
Erd\H{o}s-R\'enyi \cite{ErdosRenyi60} has given us numerous insights, surprises and
techniques and has been used to count, to establish structural
properties and to analyze algorithms.

In the standard unweighted model $G_{n,p}$, each pair of vertices
$ij$ of an $n$-vertex graph is independently declared to be an edge
with probability $p$. Equivalently, one picks a random number
$X_{ij}$ for each $ij$ in the interval $[0,1]$, i.e., a point in the
unit cube, and defines as edges all pairs for which $X_{ij} \le p$.
To get a weighted graph, we avoid the thresholding step.

In this paper, we propose the following extension to the standard
model. We have a distribution $F$ in $\R_+^N$ where $N = n(n-1)/2$
allows us a coordinate for every pair of vertices. A random point
$X$ from $F$ assigns a non-negative real number to each pair of
vertices and is thus a random weighted graph. The random graph
$G_{F,p}$ is obtained by picking a random point $X$
according to $F$ and applying a $p$-threshold to determine edges,
i.e., the edge set $E_{F,p}=\set{ij:\;X_{ij} \le p}$. It is clear that this
generalizes the standard model $G_{n,p}$ which is the special case
when $F$ is uniform over a cube.

In the special case where $F(x)=1_{x\in K}$ is the indicator function
for some convex subset $K$ of $\R_+^N$ we use the notation
$G_{K,p}$ and $E_{K,p}$. Thus to obtain $G_{K,p}$ we let $X$ be a random point in
$K$. It includes the restriction of any $L_p$ ball to the positive orthant.
The case of the simplex
\[
K=\{X \in \R^N :  \forall e, X_e \ge 0, \sum_e \a_ix_e \le L\}
\]
for some set of coefficients $\a$ appears quite interesting by itself and
we treat it in detail in Section \ref{sec:SIMPLEX}. In the weighted graph setting, it corresponds
to a random graph with a bound on the total edge weight.
In general, $F$ be could be
any distribution, but we will consider a further generalization of the cube and simplex, namely, $F$ has a logconcave density
$f$. We call this a logconcave distribution.
A function $f:\R^n \rightarrow \R_+$ is {\em logconcave} if for any
two points $x,y \in \R^n$ and any $\lambda \in [0,1]$,
\[
f(\lambda x + (1-\lambda)y) \ge f(x)^\lambda f(y)^{1-\lambda},
\]
i.e., $\ln f$ is concave.

The model appears to be
considerably more general than $G_{n,p}$. Nevertheless, can we recover
interesting general properties including threshold phenomena?

The average case analysis of algorithms for NP-hard problems was pioneered by Karp \cite{K}
and in the context of graph algorithms, the theory of random graphs has played a crucial role (see \cite{FM} for a somewhat out-dated survey).
To improve on this analysis, we need tractable distributions that provide a closer bridge between average case and worst-case.
We expect the distributions described here to be a significant platform for future research.

We end this section with a description of the model and a summary of our main results.

\subsection{The generalized model}

We consider logconcave density functions $f$ whose
support lies in the positive orthant. Let $f$ be such a density. The second moment $\sigma_e^2(f)=\E(X_{e}^2)$ along each
axis $e$
will be important. We just use $\sigma_e$ when $f$ is fixed and
simply $\sigma$ when the second moment is the same along every
axis.

Fixing only the standard deviations along the axes allows highly restricted distributions, e.g.,
the line from the origin to the vector of all $1$'s. To ensure greater "spread", we require that
the density is down-monotone, i.e., for any $x,y \in \R^N$ such that $x \le y$, we have $f(x) \ge f(y)$.
This is not a significant
restriction as the measure outside a large box will be negligible.
When $f$ corresponds to the uniform density over a convex body $K$, this means that when $x \in K$, the
box with $0$ and $x$ at opposite corners is also in $K$. It also implies that $f$ can be viewed as the
 restriction to the positive orthant of a $1$-unconditional distribution for which the density
$f(x_1, \ldots,x_N)$ stays fixed when we reflect on any subset of axes, i.e., negating subset of
coordinates keeps $f$ the same. Such distributions include, e.g.,  the $L_p$ ball for any $p$ but
also much less symmetric sets, e.g., the uniform distribution over any down-monotone convex body.
We note that sampling such distributions efficiently requires only a
function oracle, i.e., for any point $x$, we can compute a function proportional
to the density at $x$ (see e.g., \cite{LV2}).

\section{Results}
\subsection{Random graphs from logconcave densities.}
Our first result estimates the point at which $G_{F,p}$ is connected in general in terms of
$n$ and $\sigma$, the standard deviation in any direction. Our main result
is that after fixing the second moments along every axis, the threshold for connectivity can be narrowed down to within an $O(\log n)$
 factor.

\begin{theorem}\label{thm:gencon}
Let $F$ be distribution in the positive orthant with a down-monotone logconcave density. Suppose that
$\E(X_e^2)=\s_e^2$ along coordinate $e$. Let $\sm=\min\s_e$ and $\sM=\max\s_e$. Then
there exist absolute constants $c_1 < c_2$ such that
\[
\lim_{n\to\infty}\Pr(G_{F,p}\ is\ connected)=
\begin{cases}
0&p < \dfrac{c_1\sm}{n}\\ \\
1&p > \dfrac{c_2\sM\ln n}{n}
\end{cases}
\]
\end{theorem}
$F$ being so general makes this theorem quite difficult to prove. It
requires several results that are trivial in $G_{n,p}$.

The reader will notice the disparity between the upper and lower bound.
\begin{conjecture}\label{conj1}\footnote{In an early version of this paper, an abstract of which appeared in FOCS 2008, we incorrectly claimed this conjecture as a theorem.}
Let $F$ be as in Theorem \ref{thm:gencon}. Then there exists a constant $c_0$ such that if $p <c_1\sm\ln n/n$
then \whp\footnote{A sequence of events $\cE_n$ is said to occur {\em with high probability} \whp, if $\lim_{n\to\infty}\Pr(\cE_n)
\to 1$ as $n\to\infty$}\ $G_{F,p}$ has isolated vertices.
\end{conjecture}

Having proven Theorem \ref{thm:gencon} it becomes easy to prove other similar results.

\begin{theorem}\label{thm:genmatch}
Let $F$ be as in Theorem \ref{thm:gencon}. Then
there exist absolute constants $c_3 < c_4$ such that
\[
\lim_{\substack{n\to\infty\\n\ even}}\Pr(G_{F,p}\ has\ a\ perfect\ matching)=
\begin{cases}
0&p < \dfrac{c_3\sm}{n}\\ \\
1&p > \dfrac{c_4\sM\ln n}{n}
\end{cases}
\]
\end{theorem}

Finally, for this section, we mention a result on Hamilton cycles that can be
obtained quite simply from a result of Hefetz, Krivelevich and Szab\'o \cite{HKS}.
\begin{theorem}\label{tham}
Let $F$ be as in Theorem \ref{thm:gencon}. Then
there exists an absolute constant $c_6$ such that if
$$p\geq c_6\sm\frac{\ln n}{n}\cdot\frac{\ln\ln\ln n}{\ln\ln\ln\ln n}$$
then
$G_{F,p}$ is Hamiltonian \whp.
\end{theorem}
\subsection{Random Graphs from a Simplex}
We now turn to a specific class of convex bodies $K$ for which we can prove fairly tight results.
We consider the special case where $X$ is chosen
uniformly at random from the simplex
$$\S=\S_{n,L,\a}=\set{X\in \R_+^N:\;\sum_{e\in E_n}\a_eX_e\leq L}.$$
Here $N=\binom{n}{2}$ and $E_n=\binom{[n]}{2}$ and $L$ is a positive
real number and $\a_e>0$ for $e\in E_n$.

We observe first that
$G_{\S_{n,L,\a},p}$ and $G_{\S_{n,N,\a N/L},p}$ have the same distribution and so
we assume, unless otherwise stated, that $L=N$. The special case where
$\a=\bf1$ (i.e. $\a_e=1$ for $e\in E_n$) will be
easier than the general case. We will see
that in this case $G_{\S,p}$ behaves a lot like $G_{n,p}$.

Although it is
convenient to phrase our theorems under the assumption that $L=N$, we
will not always assume that $L=N$ in the main body of our proofs. It is
informative to keep the $L$ in some places, in which case we will use the notation $\S_L$ for the simplex.
In general, when discussing the simplex case, we will use $\S$ for the simplex. On the other hand, we will if necessary subscript
$\S$ by one or more of the parameters $\a,L,p$ if we need to stress their values.

We will not be able to handle completely general $\a$. We will
restrict our attention to the case where
\begin{equation}\label{aM}
\frac{1}{M}\leq \a_e\leq M\qquad\qquad for\ e\in E_n
\end{equation}
where $M=M(n)$.
An $\a$ that satisfies \eqref{aM} will be called {\em M-bounded}.

This may seem restrictive, but if we allow arbitrary $\a$ then by
choosing $E\subseteq E_n$ and making $\a_e,\,e\notin E$ very small and
$\a_e=1$ for $e\in E$ then $G_{\S,p}$ will essentially be a random
subgraph of $G=([n],E)$, perhaps with a difficult distribution.

We first discuss the connectivity threshold: We need the following
notation.
$$\a_v=\sum_{w\neq v}\a_{vw}\qquad\qquad for\ v\in [n].$$
\begin{theorem}\label{S1}\
\begin{description}
\item[(a)] Let $p=\frac{\ln n+c_n}{n}$. Then if $\a=\bf1$,
$$\lim_{n\to\infty}\Pr(G_{\S,p}\ is\
connected)=\begin{cases}0&c_n\to-\infty\\e^{-e^{-c}}&c_n\to
  c\\1&c_n\to\infty\end{cases}.$$
\item[(b)] Suppose that $\a$ is $M$-bounded and $M\leq (\ln n)^{1/4}$. Let $p_0$ be the solution to
$$\sum_{v\in [n]}\xi_v(p)=1$$
where $\xi_v(p)=\brac{1-\frac{\a_vp}{N}}^N$.
Then for any fixed $\e>0$,
$$\lim_{n\to\infty}\Pr(G_{\S,p}\ is\
connected)=\begin{cases}0&p\leq (1-\e)p\\1&p\geq(1+\e)p\end{cases}.$$
\end{description}
\end{theorem}
Our proof of part (a) of the above theorem relies on the following:
\begin{lemma}\label{lem1}
If $\a=\bf1$ and $m$ is the number of edges in $G_{\S,p}$. Then
\begin{description}
\item[(a)] Conditional on $m$, $G_{\S,p}$ is distributed as $G_{n,m}$
  i.e. it is a random graph on vertex set $[n]$ with $m$ edges.
\item[(b)] \Whp\ $m$ satisfies
$$\E(m) + \sqrt{\E(m)\om} \le m \le \E(m) + \sqrt{\E(m)\om}$$
for any $\om=\om(n)$ which tends to infinity with $n$.
\end{description}
\end{lemma}
So to prove part (a) all we have to verify is that $\E(m)\sim \frac12
n(\ln n+c_n)$ and apply known results about the connectivity
threshold for random graphs, see for example Bollob\'as
\cite{BollBook} or Janson, {\L}uczak and Ruci\'nski \cite{JLRBook}.
(We do this explicitly in Section \ref{ssec:simplexCon}).
Of course, this implies much more about $G_{\S,p}$ when $\a=\bf1$. It turns out to be
$G_{n,m}$ in disguise, where $m=m(p)$.

Our next theorem concerns the existence of a giant component i.e. one
of size linear in $n$. It is somewhat weak.
\begin{theorem}\label{Giant}
Let $\e>0$ be a small positive constant.
\begin{description}
\item[(a)] If $p\leq \frac{(1-\e)}{Mn}$ then \whp\ the maximum
  component size in $G_{\S,p}$ is $O(\ln n)$.
\item[(b)] If $p\geq \frac{(1+\e)M}{n}$ then \whp\ there is a
unique giant component in $G_{\S,p}$ of size $\geq \k n$ where $\k=\k(\e,M)$.
\end{description}
\end{theorem}
Let $\cP$ be a monotone increasing graph property. $p_0$ is a threshold for $\cP$ if
$p/p_0\to 0$ implies that $\Pr(G_{\S,p}\in\cP)\to 0$ and $p/p_0\to \infty$ implies that $\Pr(G_{\S,p}\in\cP)\to 1$.
It is an open question as to whether every monotone property has a threshold. We can make the following rather
weak statement.
\begin{theorem}\label{thresh}
If $M=O(1)$ and $\a$ is $M$-bounded then every monotone property $\cP$ has a threshold in the model $G_{\S,p}$.
\end{theorem}

We turn our attention next to the diameter of  in $G_{\S,p}$.
\begin{theorem}\label{S3}
Let $k\geq 2$ be a fixed integer. Suppose that $\a$ is $M$-bounded and for simplicity assume only that
$M=n^{o(1)}$.
Suppose that $\th$ is fixed and
satisfies $\frac{1}{k}<\th<\frac{1}{k-1}$. Suppose that
$p=\frac{1}{n^{1-\th}}$. Then \whp\
$\diam(G_{\S,p})=k$.
\end{theorem}

Our next theorem concerns spanning trees.
We say that $\a$ is {\em decomposable} if there
exist $d_v,\,v\in [n]$ such that $\a_{vw}=d_vd_w$. In which case we
define
$$d_S=\sum_{v\in S}d_v\ for\ S\subseteq V\ and\ D=d_V.$$

Let $\Lambda_X$ be weight of the minimum length spanning tree of
the complete graph $K_n$ when the edge weights are given by $X$.
\begin{theorem}\label{S2}\
If $\a$ is decomposable and $d_v\in [\om^{-1},\om],\,\om=(\ln n)^{1/10}$
for $v\in V$ and $X$ is chosen uniformly at random from $\S_{n,\a}$ then
$$\E[\Lambda_X]\sim \sum_{k=1}^\infty \frac{(k-1)!}{D^k}\sum_{\substack{S\subseteq V\\|S|=k}}\frac{\prod_{v\in S}d_v}{d_S^2}.$$
(The notation $a_n\sim b_n$ means that $\lim_{n\to\infty}(a_n/b_n)=1$,
assuming that $b_n>0$ for all $n$.)
\end{theorem}
Note that if $d_v=1$ for all $v\in [n]$ then the expression in the
theorem yields $\E[\Lambda_X]\sim \zeta(3)$.

\subsection{Random Traveling Salesman Problems}
We will also consider the use of $X$ as weights for an optimisation
problem. In particular, we will consider the Asymmetric Traveling Salesman Problem (ATSP)
in which the
weights $X:[n]^2\to\R_+$ are randomly chosen from a simplex.
We will need to make
an extra assumption about the simplex. We assume that
$$\a_{v_1,w}=\a_{v_2,w}\ for\ all\ v_1,v_2,w.$$
Under this assumption, the distribution of the weights of edges
leaving a vertex $v$ is independent of
the particular vertex $v$. We call this {\em row symmetry}.
We show that
a simple patching algorithm based on that in \cite{KS} works \whp.
\bth{tsp}
Suppose that the cost matrix $X$ of an instance of the ATSP is drawn from a row symmetric simplex
where $M\leq n^\d$, for sufficiently small $\d$.
Then there is an $O(n^3)$ algorithm that \whp\ finds a
tour that is asymptotically optimal, i.e., 
\whp\ the ratio of cost of the tour found to the optimal tour cost tends to one.
\enth
\section{Proofs: logconcave densities}\label{logprooofs}
We consider logconcave distributions restricted to the positive orthant.
We also assume they are down-monotone, i.e., if $x \ge y$ then the density function $f$ satisfies
$f(y) \ge f(x)$. We begin by collecting some well-known facts
about logconcave densities and proving some additional properties. The new properties
 will be the main tools for our subsequent analyses and allow us to deal with the
 non-independence of edges.

\subsection{Properties}
The following classical theorem summarizing basic properties of logconcave functions was proved by
Dinghas \cite{Din}, Leindler \cite{Lei} and Pr\'ekopa \cite{Pre1, Pre2}.

\begin{theorem}\label{LOGCONCAVE}
All marginals as well as the distribution function of a logconcave
function are logconcave. The convolution of two logconcave functions
is logconcave.
\end{theorem}
We will need the several results from \cite{LV1}:
A logconcave function $f:\R^m\to\R_+$ is {\em isotropic} if (i) it has mean 0 and (ii)
its co-variance matrix is the identity. It is a {\em density} if $\int_xf(x)dx=1$.
If $f$ is a density then so is $f_\l(x)=\l^mf(\l x)$. Also $\s_e(f_\l)=\s_e(f)/\l$ for all $e$.
These identities are useful for translating results on the {\em isotropic} case to a more general case.
For a function $f$ we denote its maximum value by $M_f$.
\begin{lemma}\label{15}\
\begin{description}
\item[(a)] Let $f:\R\to\R_+$ be a logconcave density function with mean $\m_f$. Then
$$ \frac{1}{8\s_f}\leq f(\mu_f)\leq M_f\leq \frac{1}{\s_f}.$$
(For a one dimensional function $f$, it is appropriate to use $\s_f=\s(f)$).
\item[(b)] Let $X$ be a random variable with a logconcave density function
$f : \R \to \R_+$.
\begin{description}
 \item[(i)]  For every $c> 0$,
\[
\Pr(f(X) \le c) \le \frac{c}{M_f}.
\]
\item[(ii)] $$\Pr(X\geq \E(X))\geq \frac{1}{e}.$$
\end{description}
\item[(c)] Let $X$ be a random point drawn from a logconcave distribution in $\R^m$. Then
$$\E(|X|^k)^{1/k}\leq 2k\E(|X|).$$
\item[(d)] If $f:\R^s\to\R_+$ is an isotropic logconcave density function then
$$M_f\geq (4e\p)^{-s/2}.$$
\end{description}
\end{lemma}
\proofend

The above lemma is from \cite{LV1}. 
Part (a) of this lemma is from Lemma 5.5. Part (bi) is Lemma 5.6(a) and Part (bii) is Lemma 5.4. Part (c) is Lemma 5.22.
Part (d) is Lemma 5.14(c).

We prove the next four lemmas with our theorems in mind.
\begin{lemma}\label{lem:basic}
Let $X$ be a random variable with a non-increasing logconcave density function $f : \R \to \R_+$.
\begin{description}
\item[(a)] For any $p \ge 0$,
\[
\Pr(x \le p) \le pM_f \le \frac{p}{\sigma_f}.
\]
\item[(b)]
For any $0\leq p\leq \s_f$,
\[
\Pr\brac{x \le p} \ge \frac{pM_f}{2} \ge \frac{p}{2\s_f}.
\]
\end{description}
\end{lemma}
\proofstart
For part (a) use $\Pr(x\leq p)=\int_{x=0}^pf(x)dx \le pM_f$ and then apply Lemma \ref{15}(a).

For part (b), we check the value of $f(p)$. If $f(p) \ge M_f/2$,  then the claim follows. If not, by Lemma \ref{15}(bii),
\[
\Pr\brac{f(X) \le \frac{M_f}{2}} \le \frac{1}{2}
\]
and so
\[
\Pr(X \le p) \ge\Pr\brac{f(X) \ge \frac{M_f}{2}}\geq  \frac{1}{2} \ge \frac{p}{2\sigma_f}
\]
as required.
\proofend

\begin{lemma}\label{lem:centroid}
Let $F:\reals^s_+\to\reals_+$ be a distribution with a down-monotone logconcave density function $f$ with support in the
positive orthant. Let $\E(X_i^2)=\sigma_i^2$
for coordinate $i$ and let $\s_\Pi=\prod_{i=1}^s\s_i$. Let $v=(v_1,\ldots,v_s)$ be the centroid of $F$.
Then $v_i \ge \sigma_i/4$ for all $i \le s$ and $f(v) \geq e^{-A_1s}/\s_\Pi$ for some absolute constant
$A_1>0$.
\end{lemma}
\proofstart
Applying Lemma \ref{15}(c) with $k=2$ gives
\[
v_i = \int_{\R_+^s} x_i f(x) \, dx \ge \frac{1}{4}\left(\int_{\R_+^s}
x_i^2 f(x)\, dx\right)^{\frac{1}{2}} \ge \frac{\sigma_i}{4}.
\]

We next prove that
\begin{equation}\label{fvf0}
f(v) \ge 2^{-2s-4}f(0).
\end{equation}
Let $H \subseteq \R^s$
be a hyperplane through $v$ that is tangent to the set
$\{x: \, f(x) \ge f(v)\}$. Let $a$ be the unit normal to $H$. The down-monotonicity of $f$ implies that $a$ is non-negative.
Let $H(t)$ denote the hyperplane parallel to $H$ at
distance $t$ from the origin. Let
\[
h(t)=\int_{H(t)}f(y)dy
\]
be the marginal of $f$ along $a$. The function $h$ is also a logconcave density and observe that $\m_h=a\cdot v$.

Consider the plane $H(a\cdot v/2)$. Let $x$ be a point on $H=H(a\cdot v)$. Since $H$ is a tangent plane $f(x) \le f(v)$. Using logconcavity,
\[
f(x/2)^2 \ge f(0)f(x)
\]
and so
\[
f(x/2) \ge \sqrt{\frac{f(0)}{f(x)}}f(x) \ge \sqrt{\frac{f(0)}{f(v)}}f(x).
\]
Therefore
\[
h(a\cdot v/2) = \int_{H(a \cdot v/2)} f(y)\, dy = \frac{1}{2^{s-1}} \sqrt{\frac{f(0)}{f(v)}}h(\m_h) \ge \frac{1}{2^{s-1}} \sqrt{\frac{f(0)}{f(v)}}\frac{1}{8\s(h)}
\]
where we have used Lemma \ref{15}(a) for the last inequality.

On the other hand, using Lemma \ref{15}(a) we have $h(a\cdot v/2) \leq M_h\leq \frac{1}{\s(h)}$ and \eqref{fvf0} follows.

Applying Lemma \ref{15}(d)
to the isotropic logconcave function
$$\hf(y_1,y_2,\ldots,y_s)=2^{-s}\s_\Pi f(|\s_1y_1|,|\s_2y_2|,\ldots,|\s_sy_s|)$$
we see that $f(0)$ which is the maximum of $f$ is
at least
$(2\pi e)^{-s/2}/\s_\Pi$. The lemma follows from (\ref{fvf0}).
\proofend

\begin{lemma}\label{lem:BIGCUB}
Let $F$ be as in Lemma \ref{lem:centroid}. Let $\sm=\min\s_i$ and $\sM=\max\s_i$.
Let $G=(V,E)$ be a random graph from $G_{F,p}$
and $S \subseteq V \times V$ with $|S|=s$. Then
\[
e^{-a_1ps/\sm} \le \Pr(S \cap E = \emptyset) \le e^{-a_2ps/\sM}
\]
where $a_1,a_2$ are some absolute constants and the lower bound requires $p < \sm/4$.
\end{lemma}
\proofstart
We consider the projection of $F$ to the subspace spanned by $S$. Let $f_S$ be the resulting density function. It is logconcave by Theorem \ref{LOGCONCAVE}.
For a point $x \in \R^s_+$, let $B(x)$ be the positive orthant at $x$, i.e.,
\[
B(x)= \{y \in \R_+^s \, : \, y \ge x\}.
\]
Let $g(x)$ be the integral of $f_S$ over $B(x)$. Then by Theorem \ref{LOGCONCAVE}, $g$ is also logconcave.
The function $h(x) = \ln g(x)$ is concave and so for $e\in S$,
\[
\frac{\partial h(x)}{\partial x_e} = \frac{\frac{\partial g(x)}{\partial x_e}}{g(x)}
\]
is nonincreasing. Therefore, it achieves its maximum at $x_e=0$, i.e.,
\[
\frac{\partial h(x)}{\partial x_e} \le \frac{\partial g(0)}{\partial x_e}
\]
since $g(0)$ = 1.
The derivative of $g$ at $x_e=0$ is simply the probability mass at $x_e=0$, i.e.,
\[
\frac{\partial g(0)}{\partial x_e}=-\int_{x_e=0}f_S(x) \, dx \le -\frac{1}{8\sM}
\]
where the inequality is from Lemma \ref{lem:basic}(a).
Thus, by concavity,
\[
h(x) \le h(0) - \frac{1}{8\sM}\sum_{e\in S} x_e
\]
and so
\[
g(x) \le e^{-\sum_{e=1}^s x_i/8\sM}.
\]
Setting $x_e = p$ for all $e\in S$, we get the first inequality of the lemma.

For the lower bound, first assume that $\sM=\sm=\s$. Let $f_S$ be the marginal of $f$ in $R^S_+$ and let  $v=(v_1,\ldots,v_s)$ be the centroid of $F_s$.
Consider the box induced by the origin and $v$. From Lemma \ref{lem:centroid},
\[
g(\s/4,\s/4,\ldots,\s/4) \ge f_S(v)(\s/4)^{s} \ge e^{-(A_1+2)s}.
\]
For $p < \s/4$, by the logconcavity of $g$ along the line from $0$ to $(\s/4,\ldots, \s/4)$,
\[
g(p,\ldots,p) \ge g(0)^{1-4p/\s} g(\s/4,\ldots,\s/4)^{4p/\s} =g(\s/4,\ldots,\s/4)^{4p/\s}\ge e^{-A_2ps/\s}.
\]
We now remove the assumption $\sM=\sm$ using scaling. Define
$$\hg(y_1,y_2,\ldots,y_s)=\s_\Pi f(\s_1y_1,\s_2y_2,\ldots,\s_sy_s).$$
$\hg$ is the density of the vector $Y$ defined by
$Y_e=X_e/\s_e$ for all $e\in S$. Thus $\E(Y_e^2)=1$ for all $e \in S$ and
$$\Pr(X_e\geq p,\,e\in S)=\Pr(Y_e\geq p/\s_e,\,e\in S)\geq \Pr(Y_e\geq p/\sm,\,e\in S)\geq e^{-A_2ps/\sm}.$$
\proofend

\begin{lemma}\label{lem:SMALLCUBE}
Let $F$ be as in Lemma \ref{lem:BIGCUB}.
Let $G=(V,E)$ be a random graph from $G_{F,p}$
and $S \subseteq V \times V$ with $|S|=s$. There exist constants $b_1 < b_2$ such that
\[
\left(\frac{b_1p}{\sM}\right)^s \le \Pr(S \subseteq E) \le \left(\frac{b_2p}{\sm}\right)^s.
\]
The lower bound requires $p \le \sm/4$.
\end{lemma}
\proofstart
We prove the lemma in the case where $\sm=\sM=\s$. The general case follows by scaling as at the end of the
proof of Lemma \ref{lem:BIGCUB}.
Consider the projection to the span of $S$ and the induced density $f_S$.
{From}  Lemma \ref{lem:centroid},
we see that for $p \le \sigma/4$,  for any point
$x$ with $0 \le x_e \le p$ for all $e\in S$, $f_S(x) \ge (4e^{A_1}\sigma)^{-s}$. The lower bound follows.

For the upper bound, assume $\sm=\sM=\s$ and project to $S$ as before.
Then consider the origin symmetric function $g$ obtained by reflecting $f$ on each axis and scaling to keep it a density, i.e.,
\[
g(x_1,\ldots,x_n) = 2^{-s}f(|x_1|,\ldots,|x_n|).
\]
This function is $1$-unconditional (i.e., reflection-invariant for the axis planes) and its covariance matrix is $\sigma^2I$.
By a theorem of Bobkov and Nazarov \cite{BobNaz03}, its maximum, $g(0) \le c^s$
for an absolute constant $c$. The bound follows.
\proofend

\subsection{Proof of Theorem \ref{thm:gencon}}\label{sec:gconn}
For a set $S$, $|S|=k$, the probability that it forms a component of $G_{F,p}$, is by Lemma \ref{lem:BIGCUB}, at most $e^{-a_2pk(n-k)/\sM}$.
Therefore,
$$\Pr(G \mbox{ is not connected}) \le \sum_{k=1}^{\lfloor n/2\rfloor} \binom{n}{k} e^{-a_2pk(n-k)/\sM}.$$
It follows that for $p \ge 3\sM\ln n/(a_2n)$, the random graph is connected \whp.

We show next that if $p\leq \sM/(3eb_2n)$ then \whp\ $|E_{F,p}|\leq n/2$ and so $G_{F,p}$ cannot be connected.
Indeed, if $p\leq \sm/(2eb_2n)$ where $b_2$ is as in Lemma \ref{lem:SMALLCUBE}) and $N=\binom{n}{2}$,
$$\Pr(|E_{F,p}|\geq n/2)\leq \binom{N}{n/2}\bfrac{b_2p}{\sm}^{n/2}\leq \frac{1}{2^{n/2}}.$$
\proofend

\subsection{Proof of Theorem \ref{thm:genmatch}}\label{sec:gmatch}
The proof of Theorem \ref{thm:gencon} shows that if $p<c_1\sm/n$ then there are isolated vertices and so we can take $c_3=c_1$.
We have no hope of getting the constants $a_1,a_2$ right here for all
$F$ and so we will be content with finding a perfect matching between
$V_1=[n/2]$ and $V_2=[n]\setminus V_1$. Applying Hall's Theorem we see
that
\begin{eqnarray*}
\Pr(G_{F,p}\ has \ no\
p.m.)&\leq&2\sum_{k=1}^{n/4}\binom{n/2}{k}\binom{n/2}{k-1}e^{-a_2k(n/2-k+1)p/\sM}\\
&\leq&2\sum_{k=1}^{n/4}\bfrac{n^2e^{2-a_2np/4\sM}}{4k^2}^k\\
&=&o(1)
\end{eqnarray*}
provided $p\geq 9\sM\ln n/(a_2n)$.
\proofend

\subsection{Proof of Theorem \ref{tham}}\label{sec:gham}
We use the following  result from \cite{HKS}: Let $G=(V,E)$ have $n$ vertices and let $d=d(n)\in [12,e^{\ln^{1/3}n}]$
be a parameter
such that with $n_0=\frac{n\ln\ln n\ln d}{\ln n\ln\ln\ln n}$:
\begin{description}
\item[P1] For every $S\subset V$, if $|S|\leq n_0/d$ then $|N(S)\geq d|S|$.\\
($N(S)$ denotes the set of vertices not in $S$ that have at least one neighbor in $S$).
\item[P2] There is an edge in $G$ between any two disjoint subsets $A,B\subset V$ such that
$|A|,|B|\geq n_0/4130$.
\end{description}
If $G$ satisfies ${\bf P_1,P_2}$ then $G$ is Hamiltonian.

So let $p=\frac{\g\sM\ln n}{n}$ where lower bounds on $\g=\g(n)$ will be exposed below.
We will use $d=\frac{\ln\ln\ln n}{\ln\ln\ln\ln n}$.
First of all,
\begin{eqnarray*}
\Pr({\bf P_1}\ fails)&\leq&\sum_{s=1}^{n_0/d}\binom{n}{s}\binom{n}{ds}e^{-a_2\g s(n-s)\ln n/n}\\
&\leq&\sum_{s=1}^{n_0/d}\brac{\frac{ne}{s}\cdot \frac{n^de^d}{d^ds^d}\cdot n^{-a_2(1-o(1))\g}}^s\\
&=&o(1)
\end{eqnarray*}
if, say, $\g\geq 2d/a_2$.

Then we have
\begin{eqnarray*}
\Pr({\bf P_2}\ fails)&\leq&\binom{n}{n_0/4130}^2e^{-a_2\g (n_0/4130)^2\ln n/n}\\
&\leq&\brac{\frac{B_1n}{n_0}\cdot n^{-\g B_2n_0/n}}^{2n_0}\qquad\qquad for\ some\ B_1,B_2>0\\
&=&o(1)
\end{eqnarray*}
if, say, $\g B_2n_0\ln n/n\geq 2\ln(n/n_0)$. This is implied by $\g\geq \frac{3\ln\ln\ln n}{B_2\ln d}$.

The theorem follows.
\section{Proofs: Simplex}\label{sec:SIMPLEX}
To apply the general results we need to compute the $\s_e$:
\beq{sigsim}
\s_e^2=\frac{\a_e}{L}\int_{x_e=0}^{L/\a_e}x_e^2\brac{1-\frac{\a_ex_e}{L}}^{N-1}dx_e=\frac{2L^2}{\a_e^2N(N+1)}.
\eeq

We can obviously do better if we take accounbt of the simpler structure of the simplex.
The following lemma represents a sharpening of Lemmas \ref{lem:BIGCUB}
and \ref{lem:SMALLCUBE} for the simplex case.
\begin{lemma}\label{leftout}\
\begin{description}
\item[(a)] If $S\subseteq E_n$ and $E_p=E(G_{\S_L,p})$,
$$\Pr(S\cap E_p=\emptyset)=\brac{1-\frac{\a(S)p}{L}}^N.$$
\item[(b)] If $S,T\subseteq E_n$ and $S\cap T=\emptyset$ and
  $|T|=o(n)$ and $\a(S)|T|p,\a(T)Np,MNp=o(L)$ then
$$\Pr(S\cap E_p=\emptyset,\,T\subseteq E_p)= (1+o(1))\brac{\prod_{e\in
    T}\a_e}\bfrac{Np}{L}^{|T|}\brac{1-\frac{\a(S)p}{L}}^N.$$
\end{description}
\end{lemma}
\proofstart\\
(a)
\begin{eqnarray}
\Pr(S\cap E_p=\emptyset)&=&\frac{\vol(\S_L\cap\set{X_e\geq p:\;e\in S})}{\vol(\S_L)}\nonumber\\
&=&\frac{(L-\a(S)p)^N/(N!\prod_{e\in E_n}\a_e)}{L^n/(N!\prod_{e\in
    E_n}\a_e)}\nonumber\\
&=&\brac{1-\frac{\a(S)p}{L}}^N.\label{Dib}
\end{eqnarray}
(b) Assume first that $S=\emptyset$. For $T'\subseteq T$ and $e\notin T'$ we have
\begin{multline*}
\Pr(e\in E_p\mid X_f,f\in
T')=1-\brac{1-\frac{\a_ep}{L-\sum_{f\in T'}\a_fX_f}}^{N-|T'|}\\ \leq
\frac{\a_e(N-|T'|)p}{L-\sum_{f\in T'}\a_fX_f}\leq\frac{\a_eNp}{L}\brac{1+\frac{2\a(T')p}{L}}.
\end{multline*}
Hence
\begin{equation}\label{UPPER}
\Pr(T\subseteq E_p)\leq \brac{\prod_{e\in
    T}\a_e}\bfrac{Np}{L}^{|T|}\exp\set{\frac{2\a(T)|T|p}{L}}.
\end{equation}
Similarly,
\begin{multline*}
\Pr(e\in E_p\mid X_f,f\in
T')\geq
\frac{\a_e(N-|T'|)p}{L-\sum_{f\in T'}\a_fX_f}\brac{1-\frac{\a_e(N-|T'|)p}{2(L-\sum_{f\in T'}\a_fX_f)}}\\
\geq
\frac{\a_eNp}{L}\brac{1-\frac{|T'|}{N}-\frac{\a_eNp}{L}}.
\end{multline*}
It follows that
\begin{equation}\label{case1}
\Pr(T\subseteq E_p)=\brac{\prod_{e\in
    T}\a_e}\bfrac{Np}{L}^{|T|}\exp\set{O\brac{\frac{|T|^2}{N}+\frac{\a(T)Np}{L}}}.
\end{equation}
Now
$$\Pr(S\cap E_p=\emptyset\mid X_e,\,e\in T)=\brac{1-\frac{\a(S)p}{L-\sum_{e\in T}\a_eX_e}}^{N-|T|}.$$
So, if $T\subseteq E_p$ then
$$\Pr(S\cap E_p=\emptyset\mid X_e,\,e\in T)\geq \brac{1-\frac{\a(S)p}{L}}^N\brac{1-\frac{2\a(S)\a(T)Np^2}{L(L-\a(T)p)}}.$$
and
$$\Pr(S\cap E_p=\emptyset\mid X_e,\,e\in T)\leq \brac{1-\frac{\a(S)p}{L}}^N\brac{1+\frac{2\a(S)|T|p}{L}}$$
Part (b) follows by combining the above two inequalities with \eqref{case1}.
\proofend
\subsection{Coupling $G_{\S,p}$ and $G_{n,m}$ when $\a=\bf1$: Proof of Lemma \ref{lem1}.}\label{cps}
The distribution $G_{\S,p}$ conditioned on
any fixed number of edges $m$ is uniform over graphs with $m$ edges
i.e. is distributed as $G_{n,m}$. This is because $\S$ is {\em
  axis-symmetric} i.e. it is invariant under permutation of
coordinates.

Let $e_{ij}$ be the indicator random variable for the event that $ij$ is an
edge of $S_{p,\bf1}$ and let $m = \sum_{i,j} e_{ij}$. Let $q=\E(e_{ij})$ so that
$\E(m)=qN$. We bound the variance of $m$.
\begin{align}
\E(m^2)-\E(m)^2 &= \sum_{ij}\E(e_{ij}^2)-\E(e_{ij})^2 + \sum_{ij
\neq kl} (\E(e_{ij}e_{kl})
- \E(e_{ij})\E(e_{kl}))\nonumber\\
&\le qN + \sum_{ij \neq kl} \Pr(X_{ij} \le p \mbox{ and } X_{kl} \leq p) - \Pr(X_{ij} \le p)\Pr(X_{kl} \le p). \label{probs}
\end{align}

It follows from Lemma \ref{leftout} that,
$$q = \Pr(X_{ij} \le p) =
1- \left(1-\frac{p}{L}\right)^N.$$
Furthermore, if $p\leq L/2$ then
\begin{multline*}
\Pr(X_{kl} \leq p \mbox{ and } X_{ij} \le p)
=1-\Pr(X_{ij}\geq p)-\Pr(X_{kl}\geq p)+\Pr(X_{ij}\geq p\ and\ X_{kl}\geq p)=\\
1-2\left(1-\frac{p}{L}\right)^N+\left(1-\frac{2p}{L}\right)^N.
\end{multline*}
Using these identities, we see that if $p\leq L/2$ then
\begin{align}
\E(m^2)-\E(m)^2 &\le qN +
\frac{N(N-1)}{2}\left(1-2\left(1-\frac{p}{L}\right)^N+\left(1-\frac{2p}{L}\right)^N
- \left(1- \left(1-\frac{p}{L}\right)^N\right)^2\right)\nonumber\\
&= qN + \frac{N(N-1)}{2}\left(\left(1-\frac{2p}{L}\right)^N -
\left(1-\frac{p}{L}\right)^{2N}\right)\nonumber\\
&\le qN.\label{ppL}
\end{align}
If $p>L/2$ then $\Pr(X_{kl} \leq p \mbox{ and } X_{ij} \le p)=1-2\left(1-\frac{p}{L}\right)^N$ and so \eqref{ppL} is still true.

Using Chebyshev's inequality,
\begin{equation}\label{Chev}
\Pr(qN + \sqrt{qN\om} \le m \le qN + \sqrt{qN\om}) =
1-o(1).
\end{equation}
This completes the proof of Lemma \ref{lem1}.
\proofend
\subsection{Connectivity for $G_{\S,p}$ when $\a=\bf1$: Proof of Theorem \ref{S1} (a)}\label{ssec:simplexCon}
Suppose first that $c_n\to c$.
Let now $L=N$ and let $p=\frac{\ln n+c_n}{n}$ and let $m=|E_p|$. Then $q$ in Section
\ref{cps} satisfies
\begin{equation}\label{pq}
p-\frac{p^2}{2}\leq q\leq p.
\end{equation}
Let
$m_0=Np-n^{2/3}$ and $m_1=Np+n^{2/3}$. Now \eqref{Chev} implies that
\whp, $m_0\leq m\leq m_1$. But then
\begin{multline*}
o(1)+e^{-e^{-c}}=o(1)+\Pr(G_{n,m_1}\ is\ connected)\leq \Pr(S_{p,\bf1}\
is\ connected)\\
\leq o(1)+\Pr(G_{n,m_2}\ is\ connected)=o(1)+e^{-e^{-c}}.
\end{multline*}
Taking limits gives the result for $c_n\to c$ and the result for
$c_n\to\pm\infty$ follows by monotonicity.
\subsection{Connectivity for $G_{\S,p}$: Proof of Theorem \ref{S1} (b)}\label{ssec:simplexConAlpha}
Applying Lemma \ref{leftout}(a) we see
that for $v,w\in [n]$,
\begin{eqnarray}
\Pr(v\ is\ isolated)&=&\xi_v(p),\label{Var1}\\ \nonumber\\
\noalign{where $\xi_v=\xi_v(p)=\brac{1-\frac{\a_vp}{N}}^N$,}\nonumber\\
\Pr(v,w\ are\
isolated)&=&\brac{1-\frac{(\a_v+\a_w-\a_{vw})p}{N}}^N\label{Var2}
\end{eqnarray}
Let $p=(1-\e)p_0$. We observe first that
\begin{equation}\label{av}
\frac{1}{2M^2}\ln n\leq \a_vp_0\leq 2M^2\ln n\qquad\qquad for\ all\ v\in [n].
\end{equation}
If the upper bound breaks for some $v\in V$, then we have
$\a_wp_0\geq 2\ln n$ and $\xi_w(p_0)\leq n^{-2}$ for all $w\in [n]$ and
this contradicts the definition of $p_0$. On the other hand, if the
lower bound breaks for some $v\in V$ then $\a_wp_0\leq \frac12\ln n$ and $\xi_w(p_0)\geq (1-o(1))n^{-1/2}$ for all
$w\in [n]$ and
this also contradicts the definition of $p_0$. It follows that
$\xi_v(p_0)=n^{-a_v}$ where
\begin{equation}\label{AV}
\frac{1}{3M^2}\leq a_v\leq 3M^2\ for\
v\in [n].
\end{equation}
Consider the function
$$\f(x)=\sum_{v\in [n]}n^{-x a_v}.$$
We know that $\f(1)=1$ and $\f'(1)=-\ln n\sum_va_vn^{-a_v}\leq
-\ln n/3M^2$. It follows that $\f(1-\e)=\OM((\ln n)^{1/2})$ for small $\e$ and
this implies that if $Z_0$ is the expected number of isolated vertices
in $G_{\S,p}$ then $\E(Z_0)=\OM((\ln n)^{1/2})$.

\para
Since $M=o(\ln n)$, \eqref{Var1} and \eqref{Var2}
imply that
$$\Pr(v,w\ are\
isolated)\sim \Pr(v\ is\ isolated)\Pr(w\ is\ isolated)$$
and then the
Chebyshev inequality implies that $Z_0\neq 0$ \whp\  and hence \whp\
$S_{n,p,\a}$ is not connected.

\para
Suppose now that $p=(1+\e)p_0$.
It follows from \eqref{AV} that the expected number of isolated
vertices $A_1$ in $G_{\S,p}$ satisfies
$$A_1=\sum_{v\in
  [n]}\xi_v(p)\leq n^{-\e/6M^2}\sum_{v\in
  [n]}\xi_v(p_0)=n^{-\e/6M^2}.$$
Thus \whp\ $G_{\S,p}$ has no isolated vertices.
Let $A_k$
denote the expected number of components of size $1\leq k\leq n/2$ in
$G_{\S,p}$. Let $\p_k=\Pr(A_k\neq 0)$ and $k_0=n/M^6(\ln n)^2$. Then
for $2\leq k\leq k_0$,
\begin{eqnarray}
\p_k&\leq&\sum_{|S|=k}\brac{1-\frac{\a(S:\bar{S})p}{N}}^N\label{conn}\\
&\leq&e^{k^2Mp}\sum_{|S|=k}\exp\set{-\sum_{v\in S}\a_vp}\nonumber\\
&\leq&\frac{e^{k^2Mp}A_1^k}{k!}\nonumber\\
&\leq&\bfrac{e^{kM(1+\e)(2M^3\ln n/n)}n^{-\e/6M^2}e}{k}^k\nonumber\\
&\leq&\bfrac{e^{1+o(1)}n^{-\e k/6M^2}}{k}^k\nonumber
\end{eqnarray}
for $k\leq k_0$, after using $p_0\leq 2M^3\ln n/n$ from \eqref{av}.
Thus $\sum_{k=1}^{k_0}A_k=o(1)$ and so \whp\ there are no components
of size $1\leq k\leq k_0$ in $G_{\S,p}$.

For $k>k_0$ we use
\begin{eqnarray*}
\sum_{k=k_0}^{n/2}\p_k&\leq&\sum_{k=k_0}^{n/2}\sum_{|S|=k}\brac{1-\frac{knp}{2MN}}^N\\
&\leq&\sum_{k=k_0}^{n/2}\binom{n}{k}e^{-k\ln n/(4M^3)}\\
&\leq&\sum_{k=k_0}^{n/2}\brac{\frac{ne}{k}\cdot n^{-1/4M^3}}^k\\
&\leq&\sum_{k=k_0}^{n/2}(M^6(\ln n)^2n^{-1/4M^3})^k\\
&=&o(1).
\end{eqnarray*}
Thus \whp\ there are no components
of size $1\leq k\leq n/2$ in $G_{\S,p}$.
This completes the proof of part (b) of Theorem \ref{S1}.
\proofend

\subsection{Giant Component in $G_{\S,p}$: Proof of Theorem
  \ref{Giant}}
We use a simple coupling argument. For a vector ${\bf
  p}\in \reals^N_+$ we define $G_{\a,{\bf p}}$ to be the random graph
where $X$ is chosen uniformly from $\S_\a$ and an edge $e$ is taken iff
$X_e\leq p_e$. Suppose first that $\l_e>0$ for all $e\in E_n$. Define
$\a'$ by $\a'_e=\a_e\l_e$ and define ${\bf p}'$ by $p_e'=p_e/\l_e$. We
claim that $G_{\a,{\bf p}}=G_{\a',{\bf p}'}$ in
distribution. Indeed, for a fixed graph $G=(V,E)$ we have
\begin{align*}
&\Pr(G_{\a,\bp}=G)\\
&\frac{1}{\vol(\S_N)}\int_{\substack{0\leq x_e\leq p_e\\e\in E}}{\bf vol}_{N-|E|}\brac{\set{
x_f\geq p_f,f\notin E,\,\sum_{f\notin E}\a_fx_f\leq N-\sum_{e\in E}\a_ex_e)}}\prod_{x\in E}dx_e\\
&=\brac{\prod_{e\in E}\a_e}\frac{N!}{(N-|E|)!L^N}\int_{\substack{0\leq x_e\leq p_e\\e\in E}}\brac{\max\set{0,N-\sum_{e\in E}\a_ex_e-\sum_{e\notin E}
\a_ep_e}}^{N-|E|}\prod_{x\in E}dx_e\\
&=\brac{\prod_{e\in E}\a_e'}\frac{N!}{(N-|E|)!L^N}\int_{\substack{0\leq y_e\leq p_e'\\e\in E}}
\brac{\max\set{0,N-\sum_{e\in E}\a_e'y_e-\sum_{e\notin E}
\a_e'p_e'}}^{N-|E|}\prod_{e\in E}dy_e\\
&=\Pr(G_{\a',\bp'}=G)
\end{align*}

So for (a) we start with ${\bf p}=p{\bf 1}$ and take $\l_e=1/\a_e$ to
get $G_{\S,p}=G_{{\bf 1},{\bf p}'}$ in distribution. Note that
$p_e'\leq (1-\e)/n$ and so we can couple so that $G_{{\bf 1},{\bf p}'}\subseteq G_{{\bf 1},\frac{1-\e}{n}{\bf 1}}$. Part (a) follows from
\eqref{Chev} as in Section \ref{cps}. Part (b) is similar.
\subsection{Thresholds: Proof of Theorem \ref{thresh}}
Let $p_0$ be defined by $\Pr(G_{\S,p_0}\in\cP)=1/2$.
We follow the strategy of the previous section and obtain $G_{\S,p}=G_{{\bf 1},{\bf p}'}$ in distribution and
\begin{equation}\label{squeeze}
G_{{\bf 1},\frac{p}{M}{\bf 1}}\subseteq G_{{\bf 1},{\bf p}'}\subseteq G_{{\bf 1},Mp{\bf 1}}
\end{equation}
Suppose now that $\om\to\infty$. Putting $p=\om p_0$ in \eqref{squeeze} we get
$$\Pr(G_{\S,\om p_0}\in\cP)\geq \Pr(G_{{\bf 1},\frac{\om p_0}{M}{\bf 1}}\in \cP)=1-o(1).$$
Putting $p=p_0/\om$ in \eqref{squeeze} we get
$$\Pr(G_{\S,\om p_0}\in\cP)\leq \Pr(G_{{\bf 1},\frac{M p_0}{\om}{\bf 1}}\in \cP)=o(1).$$
\proofend

\subsection{Diameter of $G_{\S,p}$: Proof of Theorem \ref{S2}}
Recall that $p=\frac{1}{n^{1-\th}}$ where $\frac{1}{k}<\th<\frac{1}{k-1}$.
We show first that \whp\ the diameter exceeds $k-1$. Let $Z_t$ denote
the number of paths of length $t\leq k-1$ from vertex 1 to vertex
2. We consider the existence of $t$ edges making up a path.
Applying Lemma \ref{leftout}(b): $S=\emptyset$ and $|T|=k$,
\begin{eqnarray*}
\E[Z_t]&\leq&(1+o(1))n^{t-1}(Mp)^t\\
&\leq&2n^{t-1}\bfrac{M}{n^{1-\th}}^t\\
&=&2M^tn^{\th t -1}\\
&=&o(1).
\end{eqnarray*}
{\bf Case 1:} $k\geq 3$.\\
We must now show that the diameter is at most $k$. The following lemma
provides some structure:
\begin{lemma}\label{L3}
The following hold \whp:
\begin{description}
\item[(a)] The maximum degree $\D\leq \D_0= 10Mn^{\th}$.
\item[(b)] If $S\subseteq V$ with $|S|\leq n^{1-\th-\e}$ for some fixed
  $\e$. Then $|N(S)|\geq n^\th |S|/(10M\ln n)$ where $N(S)$ is the set of
  vertices, not in $S$, that are neighbors of $S$.
\end{description}
\end{lemma}
\proofstart
(a) We consider the existence of $t=10Mn^\th$ edges incident with a fixed
vertex. Applying Lemma \ref{leftout}(b): $S=\emptyset$ and $|T|=\D_0$. ($k\geq 3$ is needed here to
ensure that $\a(T)p=o(1)$).
$$\Pr[\D\geq \D_0]\leq (1+o(1))n\binom{n}{\D_0}(Mp)^{\D_0}\leq 2n\bfrac{e}{10}^{\D_0}=o(1).$$
(b) Using Lemma \ref{leftout}(a) we see that the probability that this fails to hold can be bounded by
\begin{multline*}
\sum_{|S|=1}^{n^{1-\th-\e}}\,\sum_{|T|=0}^{n^\th
  s/(10M\ln n)}\brac{1-\frac{|S|(n-|S|-|T|)p}{MN}}^N\leq\\
\sum_{s=1}^{n^{1-\th-\e}}\,\sum_{t=0}^{n^\th
  s/(10M\ln n)}n^{s+t}\exp\set{-s(n-s-t)n^{\th-1}/M}\leq\\
\sum_{s=1}^{n^{1-\th-\e}}\,\sum_{t=0}^{n^\th
  s/(10M\ln n)}n^{s+t}e^{-sn^{\th}/2M}=o(1).
\end{multline*}
\proofend

For a vertex $v$ let $N_r(v)$ be the set of vertices at distance $r$
from $v$. Let $r_0=\rdown{\frac{k-1}{2}}$ and $r_1=\rdown{\frac{k}{2}}$. It follows from Lemma \ref{L3} that \whp\ we have
$$(n^{\th}/(10M\ln n))^r\leq |N_r(v)|\leq (10Mn^{\th})^r\qquad for\ 1\leq
r\leq r_1.$$
Furthermore, we have $r_0+r_1\leq k-1$. So suppose that
$v,w\in V$ and $N_{r_0}(v)\cap N_{r_1}(w)=\emptyset$. (If the
intersection is non-empty then their distance is already $\leq k$).
Now condition on the sets $T,S$  of edges and non-edges exposed in the
construction of $N_{r_0}(v),N_{r_1}(w)$. Then \whp\ we have
$|S|=O(n(M\D_0)^{r_1})$ and $|T|=O((M\D_0)^{r_1})$.

Let $\n_v=|N_{r_0}(v)|,\n_w=|N_{r_1}(w)|$. Given $S,T$ let $R=\set{xy:x\in N_{r_0}(v),y\in N_{r_1}(w)}$.
Using Lemma \ref{leftout}(b), the conditional probability that there is no edge between
$N_{r_0}(v)$ and $N_{r_1}(w)$ is bounded as follows: $|R|+|S|=O(n^{r_1\th+1+o(1)})$ and $|T|=O(n^{r_1\th+o(1)})$.
\begin{multline}
\frac{\Pr((R\cup S)\cap E_p=\emptyset,\,T\subseteq E_p)}{\Pr(S\cap E_p=\emptyset,\,T\subseteq E_p)}=(1+o(1))\brac{1-\frac{\a(R)p}{N}}^N\\
\leq 2e^{-\n_v\n_wp/M}=\exp\set{-\Omega(n^{(r_0+r_1+1)\th-1-o(1)})}.\label{EQ10}
\end{multline}
Now $(r_0+r_1+1)\th-1=\Omega(1)$ and this completes the proof for the case
$k\geq 3$.

{\bf Case 2:} $k=2$.\\
This is much simpler. We show that if
$p=n^{-\b}$ where $\b=1/2-\e$ then $\diam(G_{\S,p})=2$
\whp. Here $\e$ is an arbitrarily small positive constant.

We first argue that the minimum degree in $G_{\S,p}$ is at least
$\D_1=n^{1/2+\e}/(10M\ln n)$. Indeed, if $\d$ denotes minimum degree then from Lemma \ref{leftout}(a),
$$\Pr[\d\leq \D_1]\leq n\binom{n}{n-\D_1}\brac{1-\frac{(n-\D_1)p}{MN}}^N=o(1).$$
Then by conditioning on $N(v)$, we argue as in \eqref{EQ10} that
\whp\ every pair of distinct vertices $v,w$ have a common neighbour. More precisely,
$$\frac{\Pr(v,w\ have\ no\ common\ nbr, N(v)=X)}{\Pr(N(v)=X)}=(1+o(1))\brac{1-\frac{\D_1p}{MN}}^N\leq e^{-n^\e}.$$
\subsection{Minimum Spanning Tree: Proof of Theorem \ref{S3}}\label{ssec:simplexMST}
Suppose that $T$ is our minimum length spanning tree. Then we can
write its length $\ell(T)$ as
\begin{eqnarray*}
\ell(T)&=&\sum_{e\in T}X_e\\
&=&\sum_{e\in T}\int_{p=0}^N1_{X_e\geq p}dp\\
&=&\int_{p=0}^N\sum_{e\in T}|\set{e:\;X_e\geq p}|dp\\
&=&\int_{p=0}^N(\k(G_{\S,p})-1)dp
\end{eqnarray*}
where $\k$ denotes the number of components.

So,
\begin{equation}\label{EQ7}
\Lambda_X=\int_{p=0}^N(\E[\k(G_{\S,p})]-1])dp
\end{equation}
Going back to \eqref{conn} (with $M=\om^2$) we see that
\begin{equation}\label{connx}
\p_k\leq \binom{n}{k}\brac{1-\frac{knp}{2\om^2N}}^N\leq \brac{\frac{ne}{k}\cdot e^{-np/2\om^2}}^k
\end{equation}
for $1\leq k\leq n/2$.

So,
$$p\geq p_0=\frac{5\om^2\ln n}{n}\ implies\ \Pr[G_{\S,p}\ is\ not\
connected]=o(N^{-2}).$$
So,
\begin{equation}\label{EQ8}
\Lambda_X=\int_{p=0}^{p_0}(\E[\k(G_{\S,p})]-1])dp+o(N^{-1}).
\end{equation}
Next let $\k_{k,p}$ denote the number of components with $k$ vertices. $\k_{1,p}$ is the number of isolated vertices and
$$\E[\k_{1,p}]=\sum_{v\in V}\brac{1-\frac{d_v(D-d_v)p}{N}}^N.$$
It follows that
\begin{equation}\label{La0}
\Lambda_X\geq \int_{p=0}^{p_0}\sum_{v\in V}\brac{1-\frac{d_v(D-d_v)p}{N}}^Ndp-p_0+o(N^{-1})\geq \Lambda_0=\frac{1}{2D}\sum_{v\in V}\frac{1}{d_v}\geq \frac{1}{2\om^2}.
\end{equation}
Using Lemma \ref{leftout}(b) to tighten \eqref{connx}, we see that for $k\leq n^{1/2}$ and $p\leq p_0$,
\begin{equation}\label{qq}
\E[\k_{k,p}]\leq
\sum_{|S|=k}k^{k-2}(\om^2p)^{k-1}\brac{1-\frac{knp}{2\om^2N}}^N\leq
\frac{1}{\om^2p}\brac{ne\cdot \om^2pe^{-np/2\om^2}}^k.
\end{equation}
{\bf Explanation:} Choose a set $S$ of $k$ vertices and then a tree $H$ on these vertices in $k^{k-2}$ ways. $(\om^2p)^{k-1}\brac{1-\frac{kn}{2\om^2N}}^N$
bounds the probability that $H$ exists and there are no edges from $S$ to $V\setminus S$.

So if $p_1=\frac{20\om^2\ln\om}{n}$ then for $k\leq n^{1/2}$,
\begin{eqnarray*}
\int_{p=p_1}^{p_0}(\E[\k_{k,p}]-1)dp&\leq& \frac{1}{\om^2p_1}\brac{2e\om^4}^k
\int_{p=p_1}^{\infty}\brac{\frac{np}{2\om^2}e^{-np/2\om^2}}^kdp\\
&=&\frac{2}{np_1}\brac{2e\om^4}^k\int_{x=10\ln \om}^{\infty}(xe^{-x})^kdx\\
&\leq&\frac{2}{np_1}\brac{2e\om^4}^k\int_{x=10\ln \om}^{\infty}e^{-2kx/3}dx\\
&\leq&\frac{2}{np_1}\brac{2e\om^4}^k\frac{1}{k\om^{6k}}\\
&\leq&\frac{1}{\om^{k+2}}.
\end{eqnarray*}
Now for any $k$ there are fewer than $n/k$ components of size $\geq k$. So,
$$\sum_{k\geq n^{1/2}}\int_{p=p_1}^{p_0}(\E[\k_{k,p}]-1)dp\leq n^{1/2}p_0.$$
It follows from \eqref{EQ8} and \eqref{La0} that
\begin{eqnarray}
\Lambda_X&=&\int_{p=0}^{p_1}(\E[\k(G_{\S,p})]-1])dp+O\brac{\sum_{k=1}^{\infty}
\frac{1}{\om^{k+2}}+n^{1/2}p_0}+o(N^{-1})\label{first}\\
&\sim& \int_{p=0}^{p_1}\E[\k(G_{\S,p})]dp\nonumber\\
&=&\sum_{k=1}^{\om^5}\int_{p=0}^{p_1}\E[\k_{k,p}]dp+O(np_1/\om^5)\nonumber\\
&\sim&\sum_{k=1}^{\om^5}\int_{p=0}^{p_1}\E[\k_{k,p}]dp,\label{La1}
\end{eqnarray}

Now let $\t_{k,p}$ denote the number of components of $G_{\S,p}$ that
are isolated trees with $k$ vertices.
For $X\subseteq V$ we let
$A_k=\set{a\in [1,k]^k:\;\sum_{j=1}^ka_j=2k-2}$. Then, where
$q=e^{-Dp}$,
\begin{equation}\label{tk}
\E[\t_{k,p}]\sim (k-2)!p^{k-1}\sum_{a\in
  A_k}\sum_{\substack{f:[k]\to V\\f\ an\ injection}}\prod_{j=1}^k\frac{d_{f(j)}^{a_j}q^{d_{f(j)}}}{(a_j-1)!}\qquad\qquad for\ k\leq\om^5.
\end{equation}
{\bf Explanation:} We choose a degree
sequence $a_j,\,j=1,2,\ldots,k$ for our tree. Then we choose $f$ to assign vertices to the degrees.
The number of trees with this
degree sequence is $\frac{(k-2)!}{\prod_{v\in X}(a_v-1)!}$. Let $H$ be
such a tree. Going back to Lemma \ref{leftout}(b) with $T=E(H)$ and
$|S|=k(n-k)+\binom{k}{2}-k+1$ we see that the probability $H$ is an
isolated tree component is $\sim p^{k-1}\prod_{v\in
  X}d_v^{a_v}\brac{1-\frac{d_vDp}{N}}^N\sim p^{k-1}\prod_{v\in
  X}d_v^{a_v}q^{d_v}$.

We will show that the expression \eqref{tk} can be re-expressed
\begin{equation}\label{tka}
\E[\t_{k,p}]\sim
(k-2)!p^{k-1}\sum_{a\in A_k}\prod_{i=1}^k\sum_{v=1}^n\frac{d_v^{a_i}q^{d_v}}{(a_i-1)!}.
\end{equation}
Observe that the sum $\S$ on the RHS of \eqref{tka} can be expressed
$$\S=\S_1+\cdots+\S_k$$
where
$$\S_j= \sum_{a\in
  A_k}\sum_{f\in \cF_j}\psi(a,f)$$
and $\cF_j$ is the set of functions from $[k]\to V$ with a range of
size $j$ and $\psi(a,f)=\prod_{i=1}^k\frac{d_{f(i)}^{a_j}q^{d_{f(i)}}}{(a_i-1)!}$.

Thus the sum on the RHS of \eqref{tk} is equal to $\S_k$. We show next that
\begin{equation}\label{ratop}
\frac{\S_{j+1}}{\S_j}\geq n^{1-o(1)}\qquad\qquad 1\leq i<k.
\end{equation}
Observe first that
$$\frac{1}{\om^{2k}e^{k\om Dp}k!}\leq \psi(a,f)\leq \om^{2k}.$$
Our bounds $\om^{10}\leq \ln n,k\leq \om^5,p\leq p_1$ imply that $\psi(a,f)=n^{o(1)}$ for
all $a,f$. So, $\S_j=|\cF_j|n^{o(1)}=n^{j+o(1)}$. This confirms
\eqref{ratop}, which implies that $\S\sim \S_k$ and confirms
\eqref{tka}.

We re-write \eqref{tka} as
\begin{eqnarray}
\E[\t_{k,p}]&\sim&(k-2)!p^{k-1}[x^{2k-2}]\brac{\sum_{v=1}^n\sum_{r=1}^\infty\frac{q^{d_v}d_v^r}{(r-1)!}x^r}^k\nonumber\\
&=&(k-2)!p^{k-1}[x^{k-2}]\brac{\sum_{v=1}^nq^{d_v}d_ve^{d_vx}}^k\nonumber\\
&=&(k-2)!p^{k-1}\sum_{\substack{S\subseteq V\\|S|=k}}q^{d_S}\frac{d_S^{k-2}}{(k-2)!}\prod_{v\in S}d_v\label{TREES}
\end{eqnarray}
where $d_S=\sum_{v\in S}d_v$.

So,
\begin{eqnarray}
\sum_{k=1}^{\om^5}\int_{p=0}^{p_1}\E[\t_{k,p}]dp&\sim&{\sum_{k=1}^{\om^5}\sum_{\substack{S\subseteq
    V\\|S|=k}}d_S^{k-2}\prod_{v\in
  S}d_v}\int_{p=0}^{p_1} p^{k-1}e^{-d_SDp}dp\nonumber\\
&=&\sum_{k=1}^{\om^5}\sum_{\substack{S\subseteq
    V\\|S|=k}}\frac{\prod_{v\in
  S}d_v}{d_S^2D^k}\int_{x=0}^{d_SDp_1} x^{k-1}e^{-x}dx\label{LL1}\\
&\sim&\sum_{k=1}^{\om^5}\sum_{\substack{S\subseteq
    V\\|S|=k}}\frac{\prod_{v\in
  S}d_v}{d_S^2D^k}\int_{x=0}^{\infty} x^{k-1}e^{-x}dx\label{LL2}\\
&=&\sum_{k=1}^{\om^5}\frac{(k-1)!}{D^k}\sum_{\substack{S\subseteq V\\|S|=k}}
\frac{\prod_{v\in S}d_v}{d_S^2}\label{LL3}\\
&\sim&\sum_{k=1}^{\infty}\frac{(k-1)!}{D^k}\sum_{\substack{S\subseteq
    V\\|S|=k}}\frac{\prod_{v\in
  S}d_v}{d_S^2}\label{LL4}
\end{eqnarray}
{\bf \eqref{LL1} to \eqref{LL2}}: $d_SDp_1\geq 20k\ln \om$ and $x\geq
20k\ln \om$ implies that $x^{k-1}\leq e^{x/2}$. Hence
$$\int_{x=d_SDp_1}^\infty x^{k-1}e^{-x}dx\leq
\int_{x=20k\ln\om}^\infty e^{-x/2}dx=2\om^{-10k}. $$
{\bf \eqref{LL3} to \eqref{LL4}}:
$$\sum_{k=\om^5}^{\infty}\frac{(k-1)!}{D^k}\sum_{\substack{S\subseteq V\\|S|=k}}
\frac{\prod_{v\in S}d_v}{d_S^2}\leq
\sum_{k=\om^5}^{\infty}\frac{(k-1)!\om^2}{k^2D^k}\sum_{\substack{S\subseteq
    V\\|S|=k}}\prod_{v\in S}d_v\leq
\sum_{k=\om^5}^{\infty}\frac{\om^2}{k^3}=O(\om^{-13})$$
which must be compared with \eqref{La0}.

It only remains to show that if $\s_{k,p}=\k_{k,p}-\t_{k,p}$ then
\begin{equation}\label{toshow}
\sum_{k=1}^{\om^5}\int_{p=0}^{p_1}\E[\s_{k,p}]dp=o(\om^{-2}).
\end{equation}
But, arguing as in \eqref{qq} we see that for $k\leq n/2$,
$$\E[\s_{k,p}]\leq
\sum_{|S|=k}k^{k}(\om^2p)^{k}\brac{1-\frac{knp}{2\om^2N}}^N\leq
\brac{nek\cdot \om^2pe^{-np/2\om^2}}^k.$$
Hence,
$$\sum_{k=1}^{\om^5}\int_{p=0}^{p_1}\E[\s_{k,p}]dp\leq
\sum_{k=1}^{\om^5}(2ek\om^4)^k\int_{p=0}^{p_1}\brac{\frac{np}{2\om^2}e^{-np/2\om^2}}^kdp
\leq \sum_{k=1}^{\om^5}(2ek\om^4)^kp_1=n^{o(1)-1}$$
and \eqref{toshow} follows.
\proofend

\section{TSP algorithm: Proof of Theorem \ref{tsp}}
A digraph is a set of edges $(i,j)$ and these can equally well be
viewed as the set of edges of a bipartite graph. So we consider there
to be a {\em digraph view} and a {\em bipartite view}.
The algorithm consists of the following:
\begin{description}
\item[Step 1] Solve the assignment problem with cost matrix $X$
  i.e. find a minimum cost perfect matching in the bipartite view.
The
  edges $(i,\ba(i))$ of the optimal assignment form a set of vertex
  disjoint cycles $C_1,C_2,\ldots,C_k$ in the digraph view.
\item[Step 2] Assume that $|C_1|\geq |C_2|\geq \cdots\geq |C_k|$.\\
For $i=k$ down to $2$: $C_1\leftarrow C_1\oplus C_i$. ({\em Patch
  $C_i$ into $C_1$)}.

\vspace{.05in}
Here $C_1\oplus C_i$ is obtained by removing an edge $(a,b)$ from
$C_1$ and an edge $(c,d)$ from $C_i$ and adding edges $(a,d),(c,b)$ to
make one cycle. These two edges are chosen to minimise the cost
$X_{ad}+X_{cb}$.
\end{description}
Each patch reduces the number of cycles by one and so the procedure
ends with a tour.

{\bf Analysis:} \\
\begin{description}
\item[(a)] The row symmetry assumption implies that the matching
  found in Step 1 is uniformly random and so in the digraph view it
  has $O(\ln n)$ cycles \whp. We prove this as follows:
For any two permutations $\p_1,\p_2$ we have
$$\Pr(\ba(X)=\p_1)=\Pr(\ba(\p_1\p_2^{-1}X)=\p_1)=\Pr(\ba(X)=\p_2).$$
It follows that \whp\
$|C_1|=\Omega(n/\ln n)$.
\item[(b)] We next put a high probability bound on the length of the longest edge in the
solution to Step 1. There are several steps:
\begin{description}
\item[(1)] We let $\om=KM(\ln n)^2$ for some large constant $K$
and argue that \whp\ every vertex in 
$$G_{\S,p_1},\,p_1=\om/n,$$
has in-degree and out-degree at least $\om_0=L\ln n$ where $L=K^{1/2}$.

To verify the degree bounds, fix a vertex
$v$ and partition $[n]\setminus \set{v}$ into sets $V_1,\ldots,V_{\om_0}$ of
size $\sim n/\om_0$.
Using Lemma \ref{leftout}(a) we see that
$$\Pr(\exists i:\;d_{p_1}(v,V_i)=0)\leq e^{-np_1/(M\om_0)}=n^{-L}$$
where $d_p(v,V_i)$ is the number of $G_{\S,p}$ neighbors of $v$ in $V_i$.

Thus with probability at least $1-n^{-L}$, $v$ has one out-neighbor in each part of the partition.
This gives an out-degree of
at least $L\ln n$ as required. In-degree is treated similarly.
If $L\geq 2$ then the failure probability is sufficient to give the result for all $v$.
\item[(2)] We use Lemma \ref{leftout}(b) and a simple first moment
  argument to argue that if in the bipartite view we have two sets
  $S,T$ contained in different sides of the partition and
  $|S|\leq n^{2/3}$ and $|T|\leq L|S|\ln n/4$ then \whp\ the induced bipartite sub-graph on
  $S\cup T$ contains at most $L|S|\ln n/2$ edges of length $\leq p_1$.
  Indeed, if $\cB$ is the event that there are $S,T$ with more edges, then
\begin{eqnarray*}
\Pr(\cB)&\leq& (1+o(1))\sum_{s=1}^{n^{2/3}}\sum_{t=1}^{Ls\ln n/4}
\binom{n}{s}\binom{n}{t}\binom{st}{Ls\ln n/2}
\bfrac{KM^2(\ln n)^2}{n}^{Ls\ln n/2}\\
&\leq&2n\sum_{s=1}^{n^{2/3}}\bfrac{ne}{s}^s\bfrac{4en}{Ls\ln n}^{Ls\ln n/4}
\bfrac{KM^2e(\ln n)^2s}{2n}^{Ls\ln n/2}\nonumber\\
&=&2n\sum_{s=1}^{n^{2/3}}\brac{\frac{ne}{s}\cdot
\bfrac{M^4L^3e^3(\ln n)^3s}{n}^{L\ln n/4}}^s\nonumber\\
&=&o(1).
\end{eqnarray*}
\item[(3)] Now suppose that the optimum solution to Step 1 contains an
  edge $(x,y)$ of length greater than $2Mn^{-1/2}$. We grow alternating
paths from $x,y$ in a breadth first manner using edges of length $\leq p_1$. Using (b1) and (b2)
we see that the levels grow at a rate $\geq L\ln n/5$ until they
are of size at least $n^{3/5}$ say.  This will happen regardless of the
matching \ba\ produced by Step 1.
Indeed, let $S_0=\set{x}$ and in general, let $S_{i+1}=\ba^{-1}(N_p(S_i)
\setminus S_0\cup\cdots\cup S_i$.
$N_p(S)$ denotes the neighbors in $G_{F,p_1}$ of a set $S$ contained in one side of the partition.
It follows from (b1) and (b2) that $|S_{i+1}|\geq L|S_i|\ln n/5$,
as long as $|S_i|\leq n^{2/3}$. So \whp\
there exists $i_0$ such that $|S_{i_0}|\geq n^{3/5}$. Similarly, if $T_0=\set{y}$ and $T_{j+1}=
\ba(N_p(T_j))\setminus T_0\cup\cdots\cup T_j$ then \whp\ there exists $j_0$ such that
$|T_{j_0}|\geq n^{3/5}$.

We can then use Lemma \ref{leftout}(a)
to argue that \whp\ there is an edge of length at most $Mn^{-1/2}$
joining the final two levels $S,T$. Indeed
$$\Pr(\exists |S|,|T|\geq n^{3/5}:there\ is\ no\ S,T\ edge\ of\ length\ \leq Mn^{1/2})\leq
\binom{n}{n^{3/5}}^2e^{-n^{7/10}}=o(1).
$$

Then exchanging along the alternating
path adds edges of total cost at most $M n^{-1/2}+o(p_1\ln n)\leq
2M n^{-1/2}$ and removes an edge of length strictly greater than this, a
contradiction.
\end{description}
\item[(b)] It follows from the above that we can \whp\ "ignore" the edges
  of length greater than 
$$p_2=M n^{-1/4}$$ 
in our construction in Step 1. Let the edges
  of length $\leq p_2$ be denoted $E_1$ and the edges of length in the
  range $[p_2,2p_2]$ be denoted $E_2$.
We observe next that \whp\ $|E_1|\leq 10M^2n^{7/4}$. Indeed, applying \eqref{UPPER} we
see that if $t=10M^2n^{7/4}$ then
$$\Pr(|E_1|\geq t)\leq \binom{N}{t}M^t\bfrac{M}{n^{1/4}}^t\exp\set{\frac{2M^3t^2}{Nn^{1/4}}}
\leq \brac{\frac{Ne}{t}\cdot\frac{M^2}{n^{1/4}}\cdot\exp\set{\frac{2M^3t}{Nn^{1/4}}}}^t=o(1).$$

Let us now condition on the exact
  lengths of the edges in $E_1$. The distribution of remaining edges can now \whp\ be written as
$X_e'=p_2+Y_e'$ where $Y'$ is chosen uniformly from a simplex
$\S'$ in at least $N'\geq N-10M^2n^{7/4}$ dimensions and with RHS $L'\geq N-10M^3n^{7/4}-Np_2$.
\begin{description}
\item[(1)] We can now argue very simply: Choose for each $2\leq i\leq
  k$ an edge $(a_i,b_i)$ of cycle $C_i$. (If $|C_i|=1$ then
  $a_i=b_i$). Then divide $C_1$ into $k$ paths $P_1,\ldots,P_k$ of length $\sim |C_1|/k$.
Arguing as in (a1) we can show that \whp\
\beq{degree}
\mbox{each $a_i$ has at least $n_0=n^{3/4}/(2(\ln n)^3)$
$E_1\cup E_2$ out-neighbors $Q_i$ in $P_i$.}
\eeq
Indeed, fix $i$ and divide $P_i$ into
$|P_i|/(2n^{1/4}\ln n)\geq n^{3/4}/(2(\ln n)^3)$ disjoint pieces, each of size $\geq 2n^{1/4}\ln n$.
The (conditional) probability that there is no $(E_1\cup E_2)$-edge from
$a_i$ to any one of these pieces is at most $e^{-2\ln n}=n^{-2}$.
This follows by applying Lemma \ref{leftout}(a) to $\S'$.

Thus \eqref{degree} holds \whp. Now further condition on the lengths of the $E_{2}$-edges from the $a_i$ to $C_1$.
The lengths of the unconditioned edges are now determined by the uniform selection from a simplex $\S"$ with
$\sim N$ coordinates and $RHS\sim N$. Let $R_i$ be the in-neighbors of the $Q_i$ on $C_1$.
Applying Lemma \ref{leftout}(a) once more, we see that
$$\Pr(\exists i:\;there\ is\ no\ R_i:b_i\ edge)\leq (\ln n)e^{-n_0p_2/M}=o(1).$$

\item[(2)] In summary, \whp\ the cost of the patching is $O(p_2\ln
  n)=o(1/M)$. Finally, the cost of the minimum tour is $\Omega(1/M)$ \whp.
We can for example show that if we only consider edges of
length at most $\e/(Mn)$ for small constant $\e$ then \whp\ at least half of the
vertices have out-degree zero. Lemma \ref{leftout}(a) shows that the expected number of isolated
vertices is $\Omega(n)$. We can then use the Chebyshev inequality to argue that there $\Omega(n)$
isolated vertices \whp.
\end{description}
\end{description}

\section{Discussion}

Our work raises several open questions.
\begin{description}
\item[0. Connectivity Threshold.] Is $\ln n/n$ the threshold for connectivity? E.g. prove Conjecture \ref{conj1}.
\item[1. Random graphs with prescribed structure.]
We can generate interesting classes of random graphs
with prescribed structure. For example, let us consider $H$-free
subgraphs of a fixed graph $G$. Let $P_H\subseteq [0,1]^{E(G)}$ be
defined as follows: Let $H_1,H_2,\ldots,H_s$ be an enumeration of the
copies of $H$ in $G$. Fix some $p_0$. $P_H$ is the set of solutions
to a linear program.
\begin{align*}
&\sum_{e\in E(H_i)}X_e > |E(H)|p_0 \quad for\ i=1,2,\ldots,s. \\
&0 \le X_{e} \le 1, \quad \forall e\in E(G).
\end{align*}
It is easy to see that $G_{P_H,p_0}$ is $H$-free and it would be interesting
to analyze important properties of $G_{P_H,p_0}$. For example, when $H$ is the
list of all triangles of the complete graph, we get triangle-free graphs.
Similarly when $H$ is a path of length $2$, we get matchings (and we can get
matchings of any fixed graph by including only the edges as coordinates).
\ignore{
We can get graphs with a fixed degree sequence with the following convex polytope. Let $d(v)$ be the degree of $v$.
\begin{align*}
&pd(v) -\eps\le \sum_{e\in N(v)}X_e \le pd(v) \quad \forall v \in V \\
&p-\eps \le X_{e} \le 1 \quad \forall e\in E(G).
\end{align*}
}

A related question is whether this formulation can be used to generate such
$H$-free graphs uniformly at random. Logconcave distributions can be sampled,
but the thresholding process might give a (slightly?) nonuniform distribution.
\item[2. Thresholds for monotone properties] Do monotone graph properties have
sharp thresholds for logconcave densities as they do for Erd\H{o}s-R\'{e}nyi random graphs?
\item[3. Giant Component.] When does $G_{F,p}$ have a giant component? We have barely scratched the
surface of this problem.
\item[4. Smoothed Analysis.] Smoothed Analysis as proposed by Spielman and Teng \cite{ST} can
be viewed as choosing the costs $X$ uniformly from a unit ball. This is a special case of what we are
proposing and it is natural to ask what can be proved about this generalisation, e.g. for Linear Programming.
\item[5. Hamilton Cycles.] Can we remove the $\frac{\ln\ln\ln n}{\ln\ln\ln\ln n}$ factor from the proof
of Theorem \ref{tham}?
\item[6. Degree Sequence.] This is a fundamental parameter and we know very little about it.
\end{description}

\end{document}